\documentclass[twocolumn, 10pt]{article}

\usepackage{amsmath, amssymb, pgf}
\RequirePackage{fancybox}
\usepackage{tikz, pgfplots}
\usepackage{algorithmic} 
\usepackage[ruled]{algorithm}

\usepackage[numbered, framed, bw]{mcode} 

\newtheorem{satz}{Theorem}[section]
\newtheorem{defi}[satz]{Definition}\sl
\newtheorem{remark}[satz]{Remark}\rm
\newtheorem{lemma}[satz]{Lemma}

\newtheorem{cor}[satz]{Corollary}

\newcommand{\eprf}{\nopagebreak\hspace*{1em} \hfill
$\Box$\noindent\vspace{2ex}\\\phantom{}}

\renewcommand{\Re}{\operatorname{Re}}

\newcommand{\rank}{\operatorname{rank}}
\newcommand{\tr}{\operatorname{trace}}

\begin{document}

\title{Computing the stochastic $H^\infty$-norm}
\author{
        Tobias~Damm, Peter~Benner, and Jan Hauth
\thanks{T. Damm  is with University of Kaiserslautern,
Department of Mathematics, 67663 Kaiserslautern, Germany, email:
damm@mathematik.uni-kl.de}
\thanks{P. Benner is with the Max Planck Institute for Dynamics of Complex Technical Systems, Sandtorstr.\ 1, 39106 Magdeburg, Germany e-mail: benner@mpi-magdeburg.mpg.de.}
\thanks{J. Hauth is with the Fraunhofer Institute for
  Industrial Mathematics (ITWM), 67663 Kaiserslautern, Germany,
email: jan.hauth@itwm.fraunhofer.de} }
\date{March 2017}

\maketitle
\begin{abstract}
The stochastic $H^\infty$-norm is defined as the $L^2$-induced norm of
the input-output operator of a stochastic linear system. Like the
deterministic $H^\infty$-norm it is characterised by a version of the
bounded real lemma, but without a frequency domain description or a
Hamiltonian condition. Therefore, we base its computation on a
parametrised algebraic Riccati-type matrix equation. 
\end{abstract}
\section{Introduction}
The $H^\infty$-norm is a fundamental concept for asymptotically stable
deterministic linear
time invariant systems. It is equal to the input/output norm of a
system both in the frequency and the time domain. It is used in robustness analysis and serves as a
performance index in $H^\infty$ control. In model order reduction, it
is an important measure for the quality of the approximation.  There
are very efficient algorithms for the  computation of the
$H^\infty$-norm, which are based on a Hamiltonian characterization.
The most widely used among these was described in \cite{BoydBala90, BruiStei90}, but recent
progress has been made e.g.\ in \cite{GuglGuer13, BenV14, FreiSpen14, AliyBenn16}.

A stochastic version of the $H^\infty$-norm was introduced by
Hinrichsen and Pritchard in \cite{HinrPrit98}.  It has a similar range
of applications as its deterministic $H^\infty$ counterpart, but its numerical
computation has hardly been considered in the literature. A major
obstacle in transferring ideas and algorithms from the deterministic case is the lack
of a suitable frequency domain interpretation or a Hamiltonian
characterization in the stochastic setup. 

In this note we present an algorithm to compute the stochastic
$H^\infty$-norm, based on a Riccati characterization. According to
the stochastic bounded real lemma, \cite{HinrPrit98}, the norm is
given as the infimum of all $\gamma>0$ for which a given parametrized
Riccati equation has a stabilizing solution.  We check the solvability of the
Riccati equation by a Newton iteration.

The paper is structured as follows. In Section 2 we introduce
stochastic systems, define the stochastic $H^\infty$-norm and state
the stochastic bounded real lemma. We also provide a new version of
the non-strict bounded real lemma and give some new bounds for the
stabilizing solution, which are proven in appendix \ref{sec:proofs}.
In  Section 3 we describe our basic algorithm and discuss ways to make
all the steps fast.  
In Section 4 we report on numerical experiments. In particular, we
compare our algorithm with an LMI solver.
To keep the notational
burden low, we confine ourselves to the case, where only one
multiplicative noise term affects the state vector. Our results can
easily be extended to more general situations which we hint at  in appendix \ref{sec:general}.

\section{The stochastic $H^\infty$-norm}
We consider stochastic linear systems of the form
\begin{align}\label{eq:stoch_basic}
  dx&=(Ax+Bu)\,dt+Nx\,dw\;,\quad y=Cx+Du\;,
\end{align}
where $A, N\in\mathbb{R}^{n\times n}$, $B\in\mathbb{R}^{n\times m}$,
$C\in\mathbb{R}^{p\times n}$, $D\in\mathbb{R}^{p\times m}$, and 
$w=(w(t))_{t\in\mathbb{R}_+}$ is a zero mean real 
Wiener process on a 
probability space $(\Omega ,{\cal F} ,\mu )$ with respect to an
increasing  family 
$({\cal F}_t)_{t\in \mathbb{R}_+}$ of $\sigma$-algebras ${\cal F}_t
\subset{\cal F} $ (e.g.~\cite{Arno74, Oeks98}).\\
Let $L^2_w(\mathbb{R}_+,\mathbb{R}^q)$ denote the corresponding space of
non-anti\-cipating  stochastic processes $v$ with values in $\mathbb{R}^q$
and norm
\[
\|v(\cdot)\|^2_{L^2_w}:={\cal E}\left(\int_0^\infty\|v(t)\|^2dt\right)<\infty,
\]
where ${\cal E}$ denotes expectation. For initial data $x(0)=x_0$ and
input $u\in L^2_w(\mathbb{R}_+,\mathbb{R}^m)$ we denote the 
solution and the output of \eqref{eq:stoch_basic} by $x(t,x_0,u)$ and
$y(t,x_0,u)$, respectively.
\begin{defi}
  System \eqref{eq:stoch_basic} is called \emph{asymptotically mean-square-stable}, if
  \begin{align*}
  \mathcal{E}(\|x(t,x_0,0)\|^2)\stackrel{t\to\infty}\longrightarrow0\;,
  \end{align*}
  for all initial conditions $x_0$. In this case, for simplicity, we
  also call the pair $(A,N)$ asymptotically mean-square stable. 
\end{defi}
If
  $(A,N)$ is asymptotically mean-square stable, then \eqref{eq:stoch_basic}
  defines an input-output operator $\mathbb{L}:u\mapsto y$ from
  $L^2_w(\mathbb{R}_+,\mathbb{R}^m)$ to
  $L^2_w(\mathbb{R}_+,\mathbb{R}^p)$ via $u\mapsto y(\cdot,0,u)$, see \cite{HinrPrit98}. By
  $\|\mathbb{L}\|$ we denote the induced operator norm,
  \begin{align}\label{eq:sthinf}
    \|\mathbb{L}\|&=\sup_{\|u\|_{L^2_w}=1}\|y(\cdot,0,u)\|_{L^2_w}\;,
  \end{align}
which is an
  analogue of the deterministic $H^\infty$-norm. We therefore call it
  the \emph{stochastic  $H^\infty$-norm} of system \eqref{eq:stoch_basic}.

\subsection{The stochastic bounded real lemma}
\label{sec:stoch-bound-real}

The norm \eqref{eq:sthinf} can be characterized
by the stochastic bounded real lemma. To this end, we define the
quadratic (Riccati-type)
mapping $\mathcal{R}_\gamma:\mathbb{R}^{n\times
  n}\to\mathbb{R}^{n\times n}$, which depends on the parameter
$\gamma>\|D\|_2$, by
\begin{align*}
 & \mathcal{R}_\gamma(X)=A^TX+XA+N^TXN-C^TC\\&-(B^TX-D^TC)^T(\gamma^2 I-D^TD)^{-1}(B^TX-D^TC)\;.
\end{align*}
Its Fr\'echet derivative at some $X\in\mathbb{R}^{n\times n}$ is the linear
mapping $(\mathcal{R}_\gamma)'_X:\mathbb{R}^{n\times
  n}\to\mathbb{R}^{n\times n}$ given by
\begin{align}\label{eq:derivR}
  (\mathcal{R}_\gamma)'_X(\Delta)&=A_X^T\Delta+\Delta A_X+N^T\Delta N\;,
\end{align}
where $A_X=A-B(\gamma^2I-D^TD)^{-1}(B^TX-D^TC)$.\\
Writing $\mathcal{L}_A:X\mapsto A^TX+XA$ and $\Pi_N:X\mapsto N^TXN$,
we have
\begin{align*}
  (\mathcal{R}_\gamma)'_X(\Delta)&=\mathcal{L}_{A_X}(\Delta)+\Pi_N(\Delta)\;.
\end{align*}
The pair $(A,N)$ is asymptotically mean-square stable if and only if
$\sigma(\mathcal{L}_A+\Pi_N)\subset\mathbb{C}_-=\{\lambda\in\mathbb{C}\;\big|\;\Re\lambda<0\}$,
e.g.\ \cite{Damm04}. 
\begin{satz}\cite{HinrPrit98} \label{thm:sbrl} \it
Assume that $(A,N)$ is asymptotically mean-square stable.
  For $\gamma>\|D\|_2$, the following are equivalent.
  \begin{itemize}
\item[(i)]  $\|\mathbb{L}\|<\gamma$.
\item[(ii)] There exists a negative definite solution $X<0$ to the
 linear matrix inequality
\begin{align}\label{eq:sbrlLMI}
\left[
  \begin{array}{cc}
    (\mathcal{L}_A+\Pi_N)(X)-C^TC&XB-C^TD\\B^TX-D^TC&\gamma^2 I-D^TD
  \end{array}
\right]>0.
\end{align}

\item[(iii)] There exists a negative definite solution $X<0$ to the
  strict Riccati inequality
$\mathcal{R}_\gamma(X)>0$.
\item[(iv)] There exists a solution $X\le 0$ to the Riccati equation
$\mathcal{R}_\gamma(X)=0$,
such that $\sigma( (\mathcal{R}_\gamma)'_X)\subset\mathbb{C}_-$.
  \end{itemize}
\end{satz}
\begin{remark}
  A solution of the Riccati equation
$\mathcal{R}_\gamma(X)=0$,
with $\sigma( (\mathcal{R}_\gamma)'_X)\subset\mathbb{C}_-$ is called a
\emph{stabilizing solution}. If it exists, then it is uniquely
defined, and is the largest solution of the inequality
$\mathcal{R}_\gamma(X)\ge 0$, see \cite{Damm04}. We will write
$X_+(\gamma)$ for this solution. By
Theorem \ref{thm:sbrl}, the norm $\|L\|$ is the infimum of all
$\gamma$ such that $\mathcal{R}_\gamma(X)=0$ possesses a stabilizing
solution, i.e.
\begin{align}\label{eq:HinfNormRic}
    \|\mathbb{L}\|&=\inf\left\{\gamma>\|D\|_2\;\big|\; \exists
                    X<0:{\mathcal{R}_\gamma(X)=0 \text{ and }\atop\sigma( (\mathcal{R}_\gamma)'_X)\subset\mathbb{C}_-}\right\}.
  \end{align}
\end{remark}

Under a controllability assumption we can also give a nonstrict
version of Theorem \ref{thm:sbrl} for asymptotically mean-square stable systems.
We define the controllability Gramian $P$ of system
\eqref{eq:stoch_basic} as the solution of
\begin{align}\label{eq:contrGram}
AP+PA^T+NPN^T=-BB^T\;.
\end{align}
 If
the system is stable, then $P$ is nonnegative definite, $P\ge 0$.
\begin{cor} \label{cor:sbrleq}\it
Assume that $(A,N)$ is asymptotically mean-square stable
  and $P>0$ in \eqref{eq:contrGram}.
 For $\gamma>\|D\|_2$, the following are equivalent.
  \begin{itemize}
\item[(i)] $\|\mathbb{L}\|\le \gamma$.
\item[(ii)] There exists a solution $X\le 0$ to the
 linear matrix inequality
\begin{align}\label{eq:sbrlLMIeq}
\left[
  \begin{array}{cc}
    (\mathcal{L}_A+\Pi_N)(X)-C^TC&XB-C^TD\\B^TX-D^TC&\gamma^2 I-D^TD
  \end{array}
\right]\ge 0.
\end{align}
\item[(iii)] There exists a solution $X\le0$ to the Riccati equation
$\mathcal{R}_\gamma(X)=0$.
  \end{itemize}
Moreover, if $\|\mathbb{L}\|=\gamma$, then $\mathcal{R}_\gamma(X)=0$
has a largest solution $X=X_+(\gamma)$, for which
$0\in\sigma\left((\mathcal{R}_\gamma)'_X\right)\subset\mathbb{C}_-\cup
i\mathbb{R}$. 
\end{cor}

\smallskip
This result is slightly stronger than
\cite[Proposition 9.6]{DammHinr01} or \cite[Corollary 5.3.14]{Damm04}, where it was
shown that (i) implies (iii) if $(A,B)$ is controllable. In the appendix we give a
new simplified proof, which can also be modified to obtain lower
bounds  for solutions of  \eqref{eq:sbrlLMIeq} as follows.

\subsection{Inequalities for solutions of the Riccati equation}
\label{sec:bounds-solut-ricc}

\begin{lemma}\label{lemma:bound1}\it
Assume that $(A,N)$ is asymptotically mean-square stable, and $\gamma>\|D\|_2$.
 Let $P^\dagger\ge 0$ be the Moore-Penrose inverse of $P$ given by \eqref{eq:contrGram}.
If $X\le 0$ satisfies \eqref{eq:sbrlLMIeq}, then 
\begin{align}\label{eq:BTXB-bound}
0\le \tr(-B^TXB)\le m^2\gamma^2\|B^TP^\dagger B\|_2\;.
\end{align}
\end{lemma}
 \smallskip
Note that $\tr(-B^TXB) $ is monotonically
decreasing. 
 Hence, if
\eqref{eq:BTXB-bound} is violated for some $X\le 0$ and $\tilde
X\le X$, then $\tilde{X}$ cannot be a solution of
\eqref{eq:sbrlLMIeq}. This bound is particularly easy to check. \\
Alternatively, we may compare with solutions of Riccati equations from
deterministic control. Let $\mathcal{R}^{\det}_\gamma$ denote the
counterpart of $\mathcal{R}_\gamma$ with $N=0$, i.e.\
\begin{align*}
  \mathcal{R}^{\det}_\gamma(X)&=\mathcal{R}_\gamma(X)-N^TXN\;.
\end{align*}

\begin{lemma}\label{lemma:bound2}\it
  Assume that $(A,N)$ is asymptotically mean-square stable, and
  $\gamma_1\ge\gamma>\|\mathbb{L}\|$.\\
Then the Riccati equation from the deterministic case 
\begin{align}\label{eq:detRic}
 \mathcal{R}^{\det}_{\gamma_1}(X)&=0
\end{align}
possesses a smallest solution $X_-\le 0$, and $X_-\le X$ for all
solutions $X$ of \eqref{eq:sbrlLMIeq}. 
\end{lemma}

\section{Computation of the stochastic $H^\infty$-norm}

To exploit the characterization \eqref{eq:HinfNormRic}, we need a
method to check, whether the Riccati equation
$\mathcal{R}_\gamma(X)=0$ possesses a stabilizing solution. Given
the Fr{\'e}chet derivative of $R_\gamma(X)$ displayed in \eqref{eq:derivR}, it is
natural to apply Newton's method  to solve the
stochastic algebraic Riccati equation from part (iv) of Theorem \ref{thm:sbrl}. The 
following result was proven in \cite{DammHinr01}.
  \begin{satz}\it
    Let $(A,N)$ be mean-square stable and assume that $\gamma>\|\mathbb{L}\|$. Consider the Newton iteration
    \begin{align}\label{eq:Newton}
      X_{k+1}&=X_k-{(\mathcal{R}_\gamma)_{X_k}'}^{-1}(\mathcal{R}(X_k))\;,
    \end{align}
where we assume $\sigma((\mathcal{R}_\gamma)_{X_0}')\subset\mathbb{C}_-$.
   Then the sequence $X_k$ converges to $X_+$, and for all $k\ge 1$ it
   holds that
   \begin{align}\label{eq:Xk_conditions}
\sigma((\mathcal{R}_\gamma)_{X_k}')\subset\mathbb{C}_-,\quad 
      \mathcal{R}_\gamma(X_{k})\le0, \text{ and } X_{k}\ge
     X_{k+1}\,. 
   \end{align}
  \end{satz}

Moreover, under the given
assumptions $X_0=0$ is a suitable initial guess (see appendix for a proof).
\begin{lemma}\label{lemma:0isstab}\it
  Let $(A,N)$ be mean-square stable and assume that $\gamma>\|\mathbb{L}\|$. Then
  $\sigma((\mathcal{R}_\gamma)_{0}')\subset\mathbb{C}_-$. 
\end{lemma}

For a given $\gamma>\|D\|_2$, we can check whether
   $\gamma>\|\mathbb{L}\|$ by running the Newton iteration
   \eqref{eq:Newton} starting from $X_0=0$. If all iterates are
   stabilizing, and the sequence converges with a
   given level of tolerance, then we conclude that
   $\gamma\ge\|\mathbb{L}\|$.

Conversely, if $\gamma<\|\mathbb{L}\|$, then either
$\sigma((\mathcal{R}_\gamma)_{X_k}')\not\subset\mathbb{C}_-$ for some
$k$, or the sequence $X_k$ is monotonically decreasing and unbounded.

If for some $k$ the
condition $\sigma((\mathcal{R}_\gamma)_{X_k}')\subset\mathbb{C}_-$ is
violated or the iteration takes more 
   than a fixed number of steps, then we conclude that
   $\gamma\le\|\mathbb{L}\|$. Additionally we might test the
   conditions of Lemma \ref{lemma:bound1} or Lemma \ref{lemma:bound2} in each
   step and conclude that $\gamma\le\|\mathbb{L}\|$ if one of them is not
   fulfilled. However, in all our examples only the stability
   condition was relevant.

Using bisection, we can thus compute
   $\|\mathbb{L}\|$ up to a given precision. 

\subsection{The basic algorithm}

We summarize this approach as our \emph{basic algorithm}.
    \begin{algorithm}[H]
 \caption{\em Computation of the stochastic $H^\infty$-norm}
\label{algo:BASIC}
 \begin{algorithmic}[1]
   \STATE Choose $\gamma_0<\|\mathbb L\|<\gamma_1$, $k_{\max}$,
   tol \label{alg:choosegammas}
   \REPEAT 
   \STATE Set $\gamma=\frac{\gamma_0+\gamma_1}2$, $X_0=0$ 
   \REPEAT 
    \IF{$\sigma((\mathcal{R}_\gamma)_{X_k}')\subset\mathbb{C}_-$}\label{alg:alpha}
    \STATE\label{alg:Newton_step}
    $X_{k+1}=X_k-{(\mathcal{R}_\gamma)_{X_k}'}^{-1}(\mathcal{R}(X_k))$
   \ENDIF
  \UNTIL convergence \OR $k=k_{\max}$ \OR $\sigma((\mathcal{R}_\gamma)_{X_k}')\not\subset\mathbb{C}_-$
    \IF{convergence}
    \STATE $\gamma_1=\gamma$, 
    \ELSE 
    \STATE $\gamma_0=\gamma$
    \ENDIF
\UNTIL $\gamma_1-\gamma_0<\text{tol}$
 \end{algorithmic}
\end{algorithm}
The stability test in line \ref{alg:alpha} and the solution
of the linear system in line \ref{alg:Newton_step} are central issues.
Both concern the generalized Lyapunov mapping $\mathcal{R}_{X_k}'$.
A naive implementation with general purpose eigenvalue and linear
system solvers, respectively, would result in an
overall complexity of about $\mathcal{O}(n^6)$. About the same
complexity  is required for LMI-solvers.
It is, however, well
known that standard Lyapunov equations of the form $\mathcal{L}_{A_{X_k}}(X)=Y$ can be solved in
$\mathcal{O}(n^3)$ operations, using e.g.\ the Bartels-Stewart
algorithm, \cite{BartStew72}. Exploiting this in iterative approaches,
we can bring down the complexity of Algorithm \ref{algo:BASIC}
at least to  $\mathcal{O}(n^3)$. This will be explained briefly in the
following two subsections. Moreover, we suggest a way to choose
$\gamma_0$ and $\gamma_1$ in line \ref{alg:choosegammas}.

 In the numerical experiments, we will show
that our algorithm outperforms general purpose LMI methods.

\subsection{The stability test}
\label{sec:stab-test}
The condition  $\sigma((\mathcal{R}_\gamma)_{X_k}')\subset\mathbb{C}_-$
in line  \ref{alg:alpha}
holds if and only if $\sigma(A_{X_k})\subset\mathbb{C}_-$ and
$\rho(\mathcal{L}_{A_{X_k}}^{-1}\Pi_N)<1$, where $\rho$ denotes the
spectral radius, \cite[Theorem 3.6.1]{Damm04}. Hence, we can first
check, whether  $\sigma(A_{X_k})\subset\mathbb{C}_-$ and then apply
the power method to compute the spectral radius $\rho$ of
$\mathcal{L}_{A_{X_k}}^{-1}\Pi_N$. Note that the mapping
$-\mathcal{L}_{A_{X_k}}^{-1}\Pi_N$ is nonnegative, in the sense that
it maps the cone of nonnegative definite matrices to itself, see
\cite{Damm04}. Hence, the iterative scheme
    \begin{align*}
     P_0=I,\; P_{k+1}=-\mathcal{L}_{A_{X_k}}^{-1}\Pi_N(P_k), \;\rho_k=\frac{\tr(P_kP_{k+1})}{\tr(P_kP_k)}
    \end{align*}
produces a sequence of nonnegative definite matrices $P_k$ which
generically converge to the dominant eigenvector. In the limit we have
$P_{k+1}\approx\rho P_k$, i.e.\ $\rho_k\stackrel{k\to\infty}\to\rho$.

\subsection{The generalized Lyapunov equation}
\label{sec:gener-lyap-equat}
In the Newton step in line \ref{alg:Newton_step}, 
the generalized Lyapunov equation
  \begin{align}\label{eq:LyapDelta}
    A_{X_k}^T\Delta+\Delta A_{X_k}+N^T\Delta N&=-\mathcal{R}_\gamma(X_k)
  \end{align}
has to be solved for $\Delta$ to obtain
$X_{k+1}=X_k+\Delta$. Equations of this type have been studied e.g.\ in
\cite{Damm08}. 

Note that $\Delta=\Delta^T\in\mathbb{R}^{n\times n}$ satisfies the fixed point equation
\begin{align*}
 \Delta&=-\mathcal{L}_{A_{X_k}}^{-1}(\Pi_N(\Delta)+\mathcal{R}_\gamma(X_k))\;.
\end{align*}
The condition $\sigma((\mathcal{R}_\gamma)'_{X_k})\subset\mathbb{C}_-$
implies $\rho(\mathcal{L}_{A_{X_k}}^{-1}\Pi_N)<1$, where $\rho$
denotes the spectral radius. Hence the fixed point iteration
\begin{align*}
 \Delta_{j+1}&=-\mathcal{L}_{A_{X_k}}^{-1}(\Pi_N(\Delta_j)+\mathcal{R}_\gamma(X_k))
\end{align*}
is convergent. In each step this iteration only requires the solution of a standard
Lyapunov equation. 
 at a cost at most in $\mathcal{O}(n^3)$. 
The speed of convergence can be improved by using
a Krylov subspace approach like \texttt{gmres} or \texttt{bicgstab}.
For details see \cite{Damm08}. More recently, also low-rank techniques
have been considered in \cite{BennBrei13, ShanSimo16,
  KresSirk14}.

\subsection{Choosing $\gamma_0$ and $\gamma_1$}
For the bisection it is useful to find suitable upper and lower
bounds for $\|\mathbb{L}\|$.  Let $G(s)=C(sI-A)^{-1}B+D$ be the
transfer function of the deterministic system obtained from
\eqref{eq:stoch_basic} by replacing $N$ with zero. The $H^\infty$-norm
$\|G\|_{H^{\infty}}$ equals the input-output norm of this deterministic system.
Then from Theorem \ref{thm:sbrl}
we conclude $\|G\|_{H^{\infty}}\le \|\mathbb{L}\|$, because the
  inequality \eqref{eq:sbrlLMI} for a given matrix $X<0$ implies that
  the corresponding linear matrix inequality with $N=0$ holds for the
  same $X$. Hence, if $\gamma>\|\mathbb{L}\|$, then
  $\gamma>\|G\|_{H^{\infty}}$. Therefore, we choose $\gamma_0=
    \|G\|_{H^{\infty}}$ and try $\gamma_1=2\gamma_0$.
      If the Newton iteration does not converge for $\gamma_1$, then
      we replace $\gamma_0$ by $2\gamma_0$ and repeat the previous
      step, until we have $\gamma_1>\|\mathbb{L}\|$.


\section{Numerical experiments}
The following experiments were carried out on a 2011 MacBook Air with
a 1.4 GHz Intel Core 2 Duo processor and 4 GB Memory running OS X
10.11.6 using MATLAB\textsuperscript{\textregistered} version R2016b.
\subsection{Random systems}
\label{sec:random-systems}
We first consider random data $(A,N,B,C)$ produced by $\texttt{randn}$. The matrix
$A$ is made stable by mirroring the unstable eigenvalues at
$i\mathbb{R}$. Then the spectral radius $\rho$ of $\mathcal{L}_A^{-1}\Pi_N$
is estimated as described in subsection \ref{sec:stab-test} and an
update of $N$ is obtained by multiplication with $(2\rho+1)^{-1/2}$.
Thus $(A,N)$ is guaranteed to be mean-square stable. We compute the
stochastic $H^\infty$-norm by our algorithm  and compare it with the
result obtained by the MATLAB\textsuperscript{\textregistered}-function \texttt{mincx}, see appendix
\ref{sec:usage-lmi-solver}. In all our tests, the relative difference
of the computed norms lies within the chosen tolerance level. The
computing times, however, differ significantly, see Table \ref{table:rand}. While for small dimensions $n$ the
implementation of the LMI-solver seems to be superior to our implementation,
for larger $n$ the algorithmic complexity becomes relevant. For
$n>100$ the LMI-solver is impractical. 
\noindent
\begin{table}[h!]
\label{table:rand}
\caption{Averaged computing times  (in sec) 
  for random systems.}
\centering
\begin{tabular}{r|c|c|c|c|c}
  $n$&10&20&40&80&160\\\hline
LMI&0.11s&0.99s&35.72s&2030s&-\\\hline
Alg\ref{algo:BASIC} &4.43s &7.98s& 24.43s& 156.6s&1156s
\end{tabular}
\end{table}
\subsection{A heat transfer problem}
This stochastic modification of a heat transfer problem described in
\cite{BennDamm11} was also discussed in \cite{BennDamm17}.
 On the unit square $\Omega=[0,1]^2$, the heat equation $T_t=\Delta T$
 for $T=T(t,x)$
is given with Dirichlet  condition $T=u_j$, $j=1,2,3$, on three
of the boundary edges and a stochastic Robin condition $n\cdot \nabla T=(1/2 +
\dot w)T$ on the fourth edge (where $\dot w$ stands for white noise).
We measure the average value $y(t)=\int_\Omega T(t,x)\,dx$

A standard 5-point finite difference
discretization on a $k\times k$ grid
leads to a modified Poisson matrix $A\in\mathbb{R}^{n\times n}$ with
$n=k^2$ and corresponding matrices $N\in\mathbb{R}^{n\times n}$, $B\in\mathbb{R}^{n\times 3}$
$C=\tfrac{1}{n}[1,\ldots,1]\in\mathbb{R}^{n\times n}$. The
$H^\infty$-norm of this discretization $\|\mathbb{L}_n\|$  approximates the induced
input/output norm of the partial differential equation. Table
\ref{table:heat} shows the computing times (in seconds) and the
computed norms (which coincide) for the two methods. For $k>9$ the
LMI-solver took to long to be considered.



\begin{table*}[h!]
\label{table:heat}
\caption{Computing times (in sec) and results for heat equation}

\centering
\begin{tabular}{r|c|c|c|c|c|c}
  $n$&25&36&49&64&81&100\\\hline
 {\small LMI}
&3.83s&36.92s&306.3s&1810s&8631s&-\\\hline
{\small  Alg\ref{algo:BASIC}} 
&9.56s&  13.25s&   26.29s&  73.38s&  129.1s&177.3s\\\hline
$\|\mathbb{L}_n\|$&0.4724  &  0.4694 &   0.4669  &  0.4647&    0.4628&0.4611
\end{tabular}

\medskip
\begin{tabular}{r|c|c|c|c|c|c}
  $n$&121&144&169&196&225&256\\\hline
{\small Alg\ref{algo:BASIC}} 
& 366.9s & 491.9s  & 808.5s & 1538s&2068s&3888s \\\hline
$\|\mathbb{L}_n\|$&  0.4596 &0.4583&0.4570 &0.4559&0.4549&0.4540
\end{tabular}
\end{table*}

\begin{figure}\centering
%
%
\definecolor{mycolor1}{rgb}{0.00000,0.44700,0.74100}%
\begin{tikzpicture}

\begin{axis}[%
width=6cm,
height=3cm,
scale only axis,
xmin=0,
xmax=3,
xtick={0,0.472410552902147,1,2,3},
xmajorgrids,
ymin=0,
ymax=1.2,
ylabel={$\rho\text{(}\gamma\text{)}$},
ymajorgrids,
axis background/.style={fill=white}
]
\addplot [color=mycolor1,solid,line width=1.5pt,mark=asterisk,mark options={solid},forget plot]
  table[row sep=crcr]{%
0.472410552902147	1.0\\
0.519651608192362	0.0933632516761597\\
0.566892663482577	0.0898346986774229\\
0.614133718772791	0.0881872050836661\\
0.661374774063006	0.087194885737497\\
0.708615829353221	0.0865244367931478\\
0.755856884643435	0.0860399996747396\\
0.80309793993365	0.0856739480023647\\
0.850338995223865	0.085388264331352\\
0.89758005051408	0.0851597382162846\\
0.944821105804294	0.0849733241805259\\
0.992062161094509	0.0848188139114532\\
1.03930321638472	0.084689026884654\\
1.08654427167494	0.0845787623508145\\
1.13378532696515	0.0844841608945677\\
1.18102638225537	0.0844022991184325\\
1.22826743754558	0.0843309233200474\\
1.2755084928358	0.0842682693207916\\
1.32274954812601	0.0842129374850914\\
1.36999060341623	0.0841638041157159\\
1.41723165870644	0.0841199574252218\\
1.46447271399666	0.0840806504767432\\
1.51171376928687	0.0840452660700699\\
1.55895482457709	0.0840132901826095\\
1.6061958798673	0.0839842916328083\\
1.65343693515752	0.0839579063334194\\
1.70067799044773	0.0839338249737759\\
1.74791904573794	0.0839117832938322\\
1.79516010102816	0.0838915543381861\\
1.84240115631837	0.0838729422376804\\
1.88964221160859	0.0838557771803079\\
1.9368832668988	0.0838399113159057\\
1.98412432218902	0.0838252153997953\\
2.03136537747923	0.0838115760254845\\
2.07860643276945	0.0837988933301711\\
2.12584748805966	0.083787079082175\\
2.17308854334988	0.0837760550787777\\
2.22032959864009	0.0837657517977555\\
2.26757065393031	0.0837561072573847\\
2.31481170922052	0.0837470660486289\\
2.36205276451074	0.0837385785102089\\
2.40929381980095	0.0837306000227967\\
2.45653487509117	0.0837230904029421\\
2.50377593038138	0.0837160133808545\\
2.55101698567159	0.0837093361489554\\
2.59825804096181	0.0837030289703897\\
2.64549909625202	0.0836970648385135\\
2.69274015154224	0.0836914191798731\\
2.73998120683245	0.0836860695944064\\
2.78722226212267	0.0836809956276036\\
2.83446331741288	0.083676178570191\\
};
\end{axis}
\end{tikzpicture}%
%
%
\definecolor{mycolor1}{rgb}{0.00000,0.44700,0.74100}%
\begin{tikzpicture}

\begin{axis}[%
width=6cm,
height=3cm,
scale only axis,
xmin=0,
xmax=3,
xtick={0,0.472410552902147,1,2,3},
xticklabels={{0},{$\|\mathbb{L}\|$},{1},{2},{3}},
xlabel={$\gamma$},
xmajorgrids,
ymin=-35,
ymax=5,
ylabel={$\alpha\text{(}\gamma\text{)}$},
ymajorgrids,
axis background/.style={fill=white}
]
\addplot [color=mycolor1,solid,line width=1.5pt,mark=asterisk,mark options={solid},forget plot]
  table[row sep=crcr]{%
0.472410552902147	0.00\\
0.519651608192362	-14.7768593759927\\
0.566892663482577	-19.5219281466351\\
0.614133718772791	-22.4891562419623\\
0.661374774063006	-24.5646167980922\\
0.708615829353221	-26.1035416933122\\
0.755856884643435	-27.2886509861073\\
0.80309793993365	-28.2266156346518\\
0.850338995223865	-28.9848128066981\\
0.89758005051408	-29.608193288129\\
0.944821105804294	-30.1279937333793\\
0.992062161094509	-30.5666239195254\\
1.03930321638472	-30.9405816747757\\
1.08654427167494	-31.262276236269\\
1.13378532696515	-31.5412132813705\\
1.18102638225537	-31.7847897873347\\
1.22826743754558	-31.9988415414118\\
1.2755084928358	-32.1880289573006\\
1.32274954812601	-32.3561143806688\\
1.36999060341623	-32.5061648973635\\
1.41723165870644	-32.6407029640978\\
1.46447271399666	-32.7618198414996\\
1.51171376928687	-32.8712620856933\\
1.55895482457709	-32.9704982466882\\
1.6061958798673	-33.0607708361529\\
1.65343693515752	-33.1431372024418\\
1.70067799044773	-33.2185019617606\\
1.74791904573794	-33.2876429377082\\
1.79516010102816	-33.3512320640943\\
1.84240115631837	-33.4098523465005\\
1.88964221160859	-33.4640117153081\\
1.9368832668988	-33.5141544087965\\
1.98412432218902	-33.5606703801556\\
2.03136537747923	-33.6039031132178\\
2.07860643276945	-33.6441561489758\\
2.12584748805966	-33.6816985615836\\
2.17308854334988	-33.7167695737201\\
2.22032959864009	-33.7495824632327\\
2.26757065393031	-33.7803278833354\\
2.31481170922052	-33.8091766952979\\
2.36205276451074	-33.8362823940942\\
2.40929381980095	-33.8617831927809\\
2.45653487509117	-33.8858038195967\\
2.50377593038138	-33.9084570723509\\
2.55101698567159	-33.9298451669623\\
2.59825804096181	-33.9500609108448\\
2.64549909625202	-33.9691887267556\\
2.69274015154224	-33.9873055485509\\
2.73998120683245	-34.0044816069146\\
2.78722226212267	-34.0207811202735\\
2.83446331741288	-34.0362629037946\\
};
\end{axis}
\end{tikzpicture}%
\caption{Spectral radius  and spectral abscissa of $X_+(\gamma)$ close
  to the critical value $\gamma=\|\mathbb{L}_{25}\|=0.47241$.}
\label{fig:alpharho}
\end{figure}
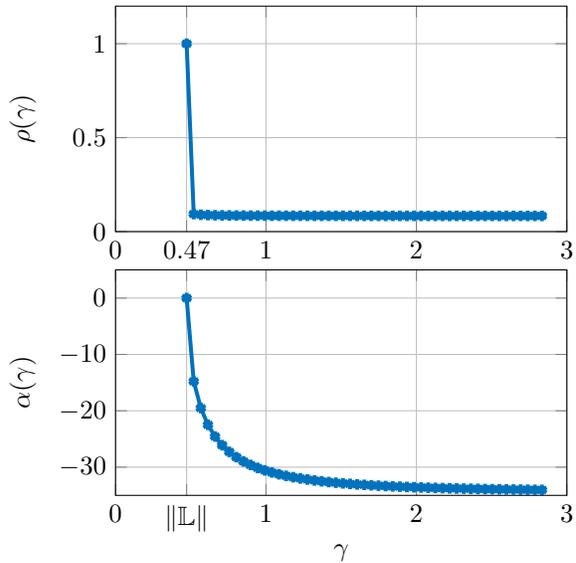

Again, we observe that Algorithm \ref{algo:BASIC} allows to treat
larger dimensions than the LMI-solver. However, the computing times
for our algorithm also grow fairly fast. As an alternative to
bisection one might consider extrapolating the spectral radii
$\rho(\gamma)=\rho\left((\mathcal{R}_\gamma)'_{X_+(\gamma)}\right)$ 
which are computed in the course of the process for $\gamma> \|\mathbb{L}\|$, or perhaps
the spectral abscissae 
$\alpha(\gamma)=\max\Re\sigma\left((\mathcal{R}_\gamma)'_{X_+(\gamma)}\right)$.
Then the norm $\|\mathbb{L}\|$ is given as the value of $\gamma$,
where $\rho(\gamma)=1$, or $\alpha(\gamma)=0$. Unfortunately, the
slopes of $\rho$ and $\alpha$ are very steep as $\gamma$ approaches
$\|\mathbb{L}\|$. Thus an extra\-polation does not seem promising. The
behaviour is visualized for the heat equation system with $n=25$ in
Figure \ref{fig:alpharho}.
\section{Conclusions}
\label{sec:conclusions}
We have suggested an algorithm to compute the stochastic
$H^\infty$-norm. 
It builds upon several ideas developed in the literature, and  is the
first algorithm, whose complexity is  considerably smaller than that of
a general purpose LMI-solver. 
We chose to present the algorithm for
the simplest case of just one multiplicative noise term, which however
can easily be generalized to the class described in appendix
\ref{sec:general}. Already in the simple
case, the stochastic $H^\infty$-norm is much harder to compute than
the $H^\infty$-norm of a deterministic system, and the computing times
are still very high. We see it as a challenge to come up with a faster
method. 

The note also contains some extensions
of known results with new
proofs, like the nonstrict stochastic bounded real lemma and lower
bounds for Riccati solutions.
\newpage
 \appendix
\section{Appendix}
\subsection{Generalization}
\label{sec:general}
System \eqref{eq:stoch_basic} can be generalized in
a straight-forward manner to the case of multiple noise terms at the
state and the input (see e.g.\ \cite{Damm04}). Then our system takes the form
\begin{align}\label{eq:stoch_general}
  dx&=(Ax+Bu)\,dt+\sum_{j=1}^\nu (N_{x,j}x+N_{u,j}u)\,dw_j\\ y&=Cx+Du\;,
\end{align}
where $N_{x,j}\in\mathbb{R}^{n\times n}$,
$N_{u,j}\in\mathbb{R}^{n\times m}$ and the $w_j$ are independent
Wiener processes.
 The Riccati operator $\mathcal{R}_\gamma$ then takes
the form
\begin{align*}
  \mathcal{R}_\gamma(X)&=P(X)-S(X)^TQ_\gamma(X)^{-1}S(X)\;,\text{
                         where }
\end{align*}
\begin{align*}
P(X)&=A^TX+XA+\sum_{j=1}^\nu N_{x,j}^TXN_{x,j}-C^TC\;,\\
S(X)&=B^TX+\sum_{j=1}^\nu N_{u,j}^TXN_{x,j}-C^TD\;,\\
Q_\gamma(X)&=\sum_{j=1}^\nu N_{u,j}^TXN_{u,j}+\gamma^2 I-D^TD\;.
\end{align*}
Our basic algorithm and all our considerations carry over to this case
literally. Only the expressions for $\mathcal{R}_\gamma$ and
$(\mathcal{R}_\gamma)'_X$ become more technical. 
 \subsection{Proofs} 
\label{sec:proofs}
As above, we write $$\mathcal{L}_A:X\mapsto A^TX+XA\quad\text{ and
}\quad\Pi_N:X\mapsto N^TXN\;.$$ On the space of symmetric matrices we
consider the scalar product
$
  \langle X,Y\rangle =\tr XY
$,
and note that the corresponding adjoint operators are 
$$\mathcal{L}_A^*:X\mapsto AX+XA^T\quad\text{ and
}\quad\Pi_N^*:X\mapsto NXN^T\;.$$ Further facts on Riccati- and Lyapunov-type
operators are cited from \cite{Damm04}.

\noindent\hspace{2em}{\itshape  Proof of Corollary \ref{cor:sbrleq}:}
(iii)$\Rightarrow$(ii) follows from the definiteness criterion via the
Schur-complement.\\
(ii)$\Rightarrow$(iii): If (ii) holds, then
$\mathcal{R}_\gamma(X)\ge0$, and by \cite{DammHinr01} there exists
 a solution $X_+\le 0$ to the equation $\mathcal{R}_\gamma(X)=0$.\\
(ii)$\Rightarrow$(i): If \eqref{eq:sbrlLMI} holds and we replace $C$ and
$D$ by $C_\varepsilon=\left[
  \begin{smallmatrix}
    C\\\varepsilon I
  \end{smallmatrix}
\right]$ and $D_\varepsilon=\left[
  \begin{smallmatrix}
    D\\0
  \end{smallmatrix}
\right]$, then we get 
\begin{align*}
\left[
  \begin{array}{cc}
    A^TX+XA+N^TXN-
C_\varepsilon^TC_\varepsilon&XB-C_\varepsilon^TD_\varepsilon\\B^TX-D_\varepsilon^TC_\varepsilon&\gamma^2 I-D_\varepsilon^TD_\varepsilon
  \end{array}
\right]> 0\;.
\end{align*}
This implies $\|\mathbb{L}_\varepsilon\|<\gamma$ for the corresponding
modified input-output operator. By
$\|\mathbb{L}_\varepsilon\|\to\|\mathbb{L}\|$ as $\varepsilon\to 0$, we
obtain $\|\mathbb{L}\|\le\gamma$.\\
(i)$\Rightarrow$(ii): If (i) holds, then $\|\mathbb{L}\|<\gamma+\frac1k$ for all
$k\in\mathbb{N}$, $k>0$. Hence there exist stabilizing solutions
$X_k\le 0$ of $\mathcal{R}_{\gamma+\frac1k}(X)=0$. Moreover, $X_k$ is
the largest solution of \eqref{eq:sbrlLMIeq} with $\gamma$ replaced by
$\gamma+\frac1k$. Hence it follows  that $X_{k+1}\le X_k$ for all $k$.
  If the $X_k$ are bounded below, then the sequence $(X_k)$ converges
  and the limit satisfies the nonstrict linear matrix inequality in
  (iii). Thus it suffices to show boundedness. We assume that the
  sequence is not bounded, i.e.\ $\|X_k\|\to\infty$ for $k\to \infty$.
  Consider the normalized sequence $\tilde X_k=\frac{X_k}{\|X_k\|}$, which -- by
  Bolzano-Weierstrass -- has a convergent subsequence $\tilde X_{k_j}$
  with limit $\tilde X\neq 0$. Then 
\begin{align*}
0&   \le 
      \frac1{\|X_{k_j}\|} 
\left[\begin{array}{cc}
(\mathcal{L}_A+\Pi_N)(X_{k_j})
        -C^TC&X_{k_j}B-C^TD\\B^TX_{k_j}-D^TC&\gamma^2 I-D^TD
      \end{array}\right]\\&
\stackrel{j\to\infty}\to 
\left[\begin{array}{cc}
        A^T\tilde X+\tilde XA+N^T\tilde X N&\tilde XB\\B^T\tilde X&0      
\end{array}\right]
\ge 0\;,
  \end{align*}
 implying $B^T\tilde X=0$ and $0\neq A^T\tilde X+\tilde XA+N^T\tilde X
N\ge 0$.   Since, by assumption $P>0$, we obtain
\begin{align*}
0&>\tr \left(P(A^T\tilde X+\tilde XA+N^T\tilde
   X N)\right)\\&=\tr \left((AP+PA^T+NPN^T)\tilde X\right)\\&=-\tr BB^T\tilde X=0
\end{align*}
which is a contradiction.\\
Thus, $\mathcal{R}_\gamma(X)=0$ has a solution $X_\infty$, which is the
limit of the largest and stabilizing solutions $X_k$ of
$\mathcal{R}_{\gamma+\frac1k}(X)\ge 0$. Thus $X_\infty$ is the largest
solution of $\mathcal{R}_\gamma(X)=0$ and $\sigma(\mathcal{R}_\gamma)'_{X_\infty} \subset\mathbb{C}_-\cup
i\mathbb{R}$.  If $\gamma=\|\mathbb{L}\|$ then $\sigma(\mathcal{R}_\gamma)'_{X_\infty} \cap
i\mathbb{R}\neq\emptyset$  and \cite[Theorem
3.2.3]{Damm04} yields that $0\in \sigma(\mathcal{R}_\gamma)'_{X_\infty}$.
\eprf

\noindent\hspace{2em}{\itshape  Proof of Lemma \ref{lemma:bound1}:}

The controllability Gramian is given by $$P=-(\mathcal{L}_A+\Pi_N)^{-*}(BB^T)\;.$$
In  the following consider an
arbitrary matrix $X\le 0$, $X\neq 0$, satisfying
$(\mathcal{L}_A+\Pi_N)(X)=Y\ge 0$. Then
 \begin{align*}
 m\|B^TXB\|_2 &\ge |\tr(B^TXB)| \\&=
\langle   (\mathcal{L}_A+\Pi_N)^{-1}(Y),-BB^T\rangle=\langle Y,P\rangle\;.
 \end{align*}
There exists a vector $u\in\mathbb{R}^m$ with $\|u\|_2=1$ and $u^*B^TXBu=-\|B^TXB\|_2\;.$
Moreover 
\begin{align}\label{eq:alphadagger}
u^TB^TYBu=\langle Y,Buu^TB^T\rangle\le \langle Y,BB^T\rangle\le \alpha_*\langle
Y,P\rangle 
\end{align}
for $\alpha_*=\|B^TP^\dagger B\|_2$. To see this,
note that the image of $B$ is contained in the image of $P$. Hence
there exists a unitary $U$, such that
\begin{align*}
  \alpha P-BB^T&=U\left[
                 \begin{array}{cc}
                   \alpha P_1-B_1B_1^T&0\\0&0
                 \end{array}
\right]U^T,\quad \text{$\det P_1\neq 0$.}
\end{align*}
The largest zero of $\chi(\alpha)=\det( \alpha
P_1-B_1B_1^T)$ is $$\alpha_*=\|P^{-1/2}B_1\|_2^2=\|B^TP^\dagger B\|_2^2\;.$$ For $\alpha\ge\alpha_*$, we have $\alpha P-BB^T\ge 0$
which proves \eqref{eq:alphadagger}.\\
We
set $\mu(X)=\langle Y,P\rangle=|\tr(B^TXB)|$.
Let now $X$ satisfy \eqref{eq:sbrlLMIeq}.
With the given data and $\eta>0$ this implies 
\begin{align*}
0&\le\left[
  \begin{array}{c}
    Bu\\\eta u
  \end{array}
\right]^*\left[
  \begin{array}{cc}
   Y&XB\\B^TX&\gamma^2 I
  \end{array}
\right]\left[
  \begin{array}{c}
    Bu\\\eta u
  \end{array}
\right]\\&
=u^*B^TYBu+2\eta u^*B^TXBu+\gamma^2\eta^2\\
&\le \alpha_*\mu(X)-\frac2m\mu(X)\eta +\gamma^2\eta^2\\
&=\gamma^2\left(\eta-\frac{\mu(X)}{
   m\gamma^2}\right)^2-\frac{\mu(X)^2}{m^2\gamma^2}+\mu(X)\alpha_*
\end{align*}
If we assume $\mu(X)> m^2\gamma^2\alpha_*$, then
the right hand is negative for 
$\eta=\frac{\mu(X)}{
   m\gamma^2}$, which is a contradiction.\\
Hence, we have $|\tr(B^TXB)|\le m^2\gamma^2 \|B^TP^\dagger B\|_2$.
\eprf
\noindent\hspace{2em}{\itshape  Proof of Lemma \ref{lemma:bound2}:}
Note that $\mathcal{R}^{\det}_{\gamma_1}(X)\ge
\mathcal{R}_{\gamma}(X)$ if $X\le 0$ and $\gamma\le\gamma_1$ and thus every solution
of $\mathcal{R}_{\gamma_1}(X)>0$ also satisfies
$\mathcal{R}^{\det}_{\gamma_1}(X)>0$. Hence
$\mathcal{R}^{\det}_{\gamma_1}(X)=0$ possesses a stabilizing solution,
and, consequently, also an anti-stabilizing solution $X_-$, which is
the smallest solution of $\mathcal{R}^{\det}_{\gamma_1}(X)\ge 0$. Thus also $X_-\le X$ for every solution $X$ of $\mathcal{R}_{\gamma_1}(X)\ge0$.
\eprf
\noindent\hspace{2em}{\itshape  Proof of Lemma \ref{lemma:0isstab}:}
We exploit the concavity of $\mathcal{R}_\gamma$ and the resolvent
positivity of $(\mathcal{R}_\gamma)'_0$, see \cite{Damm04}. If
$\|\mathbb{L}\|\le\gamma$, then there exists $X\le 0$ such that, by concavity,
\begin{align}\label{eq:concave0}
  0&= \mathcal{R}_\gamma(X)\le
     \mathcal{R}_\gamma(0)+(\mathcal{R}_\gamma)'_0(X)\;. 
\end{align}
 Assume that $\sigma
\left((\mathcal{R}_\gamma)'_0\right)\not\subset\mathbb{C}_-$. Then by
\cite[Theorem 3.2.3]{Damm04} there exists $H\ge 0$, $\lambda\ge 0$, such that
$(\mathcal{R}_\gamma)'_0(H)=\lambda H$. Taking
the scalar product of inequality \eqref{eq:concave0} with $H$, we get
\begin{align*}
  0&\le \langle \mathcal{R}_\gamma(0),H\rangle +\lambda\langle
     X,H\rangle\le 0\;.
\end{align*}
It follows that $\mathcal{R}_\gamma(0)H=0$, which implies $D^TCH=0$ and
thus $A_0H=AH$. But then $(\mathcal{L}_A+\Pi_N)^*(H)=\lambda H$ in
contradiction to the stability of $(A,N)$. \eprf

\subsection{Usage of LMI-solver}
\label{sec:usage-lmi-solver}
The LMI-solver was used as in the following
MATLAB\textsuperscript{\textregistered} listing.
\begin{lstlisting}
setlmis([])
X = lmivar(1,[n,1]);g = lmivar(1,[1,1]);
lmiterm([1 1 1 X],N',N);
lmiterm([1 1 1 X],A',1,'s');
lmiterm([1 1 1 0],C'*C);
lmiterm([1 1 2 X],1,B,'s');
lmiterm([1 2 2 g],-1,1);
lmisys = getlmis;
c = mat2dec(lmisys,zeros(n),1);
options = [tol,0,0,0,1];
copt = mincx(lmisys,c,options);
gamma = sqrt(copt)
\end{lstlisting}

\bibliographystyle{IEEEtran}
\bibliography{StochHinfNorm}

\begin{thebibliography}{10}
\providecommand{\url}[1]{#1}
\csname url@rmstyle\endcsname
\providecommand{\newblock}{\relax}
\providecommand{\bibinfo}[2]{#2}
\providecommand\BIBentrySTDinterwordspacing{\spaceskip=0pt\relax}
\providecommand\BIBentryALTinterwordstretchfactor{4}
\providecommand\BIBentryALTinterwordspacing{\spaceskip=\fontdimen2\font plus
\BIBentryALTinterwordstretchfactor\fontdimen3\font minus
  \fontdimen4\font\relax}
\providecommand\BIBforeignlanguage[2]{{%
\expandafter\ifx\csname l@#1\endcsname\relax
\typeout{** WARNING: IEEEtran.bst: No hyphenation pattern has been}%
\typeout{** loaded for the language `#1'. Using the pattern for}%
\typeout{** the default language instead.}%
\else
\language=\csname l@#1\endcsname
\fi
#2}}

\bibitem{BoydBala90}
S.~Boyd and V.~Balakrishnan, ``A regularity result for the singular values of a
  transfer matrix and a quadratically convergent algorithm for computing its
  {$L^{\infty}$}-norm,'' \emph{Syst. Control Lett.}, vol.~15, no.~1, pp. 1--7,
  1990.

\bibitem{BruiStei90}
N.~A. Bruinsma and M.~Steinbuch, ``A fast algorithm to compute the
  {$H^{\infty}$}-norm of a transfer function matrix,'' \emph{Syst. Control
  Lett.}, vol.~14, pp. 287--293, 1990.

\bibitem{GuglGuer13}
N.~Guglielmi, M.~G\"urb\"uzbalan, and M.~L. Overton, ``Fast approximation of
  the {$H_\infty$}-norm via optimization over spectral value sets,'' \emph{SIAM
  J. Matrix Anal. Appl.}, vol.~34, no.~2, pp. 709--737, 2016.

\bibitem{BenV14}
P.~Benner and M.~Voigt, ``A structured pseudospectral method for
  {$\mathcal{H}_\infty$}-norm computation of large-scale descriptor systems,''
  \emph{Math. Control Signals Syst.}, vol.~26, no.~2, pp. 303--338, 2014.

\bibitem{FreiSpen14}
M.~A. Freitag, A.~Spence, and P.~{Van Dooren}, ``Calculating the
  {$H_\infty$}-norm using the implicit determinant method,'' \emph{Linear
  Algebra Appl.}, vol.~35, no.~2, pp. 619--635, 2014.

\bibitem{AliyBenn16}
N.~Aliyev, P.~Benner, E.~Mengi, and M.~Voigt, ``Large-scale computation of
  {$\mathcal{H}_\infty$}-norms by a greedy subspace method,'' \emph{submitted
  to SIAM J. Matrix Anal. Appl.}, 2016.

\bibitem{HinrPrit98}
D.~Hinrichsen and A.~J. Pritchard, ``Stochastic ${H}_{\infty}$,'' \emph{{SIAM}
  J. Control Optim.}, vol.~36, no.~5, pp. 1504--1538, 1998.

\bibitem{Arno74}
L.~Arnold, \emph{Stochastic Differential Equations: Theory and Applications.
  Translation.}\hskip 1em plus 0.5em minus 0.4em\relax New York etc.: John
  Wiley and Sons Inc., 1974.

\bibitem{Oeks98}
B.~Oeksendal, \emph{Stochastic Differential Equations}, 5th~ed.\hskip 1em plus
  0.5em minus 0.4em\relax Springer-Verlag, 1998.

\bibitem{Damm04}
T.~Damm, \emph{Rational Matrix Equations in Stochastic Control}, ser. Lecture
  Notes in Control and Information Sciences.\hskip 1em plus 0.5em minus
  0.4em\relax Springer, 2004, no. 297.

\bibitem{DammHinr01}
T.~Damm and D.~Hinrichsen, ``Newton's method for a rational matrix equation
  occuring in stochastic control,'' \emph{Linear Algebra Appl.}, vol. 332--334,
  pp. 81--109, 2001.

\bibitem{BartStew72}
R.~H. Bartels and G.~W. Stewart, ``{A}lgorithm 432: {T}he solution of the
  matrix equation ${AX}+{XB}={C}$,'' vol. 15(9), pp. 820--826, 1972.

\bibitem{Damm08}
T.~Damm, ``Direct methods and {ADI}-preconditioned {K}rylov subspace methods
  for generalized {L}yapunov equations,'' \emph{Numer. Lin. Alg. Appl.},
  vol.~15, no.~9, pp. 853--871, 2008.

\bibitem{BennBrei13}
P.~{Benner} and T.~{Breiten}, ``Low rank methods for a class of generalized
  {L}yapunov equations and related issues.'' \emph{Numer. Math.}, vol. 124,
  no.~3, pp. 441--470, 2013.

\bibitem{ShanSimo16}
S.~Shank, V.~Simoncini, and D.~Szyld, ``Efficient low-rank solutions of
  generalized {L}yapunov equations,'' \emph{Numer. Math.}, vol. 134, pp.
  327--342, 2016.

\bibitem{KresSirk14}
D.~Kressner and P.~Sirkovi\'c, ``Greedy low-rank methods for solving general
  linear matrix equations,'' ANCHP, MATHICSE, EPF Lausanne, Switzerland,''
  Technical report, 2014.

\bibitem{BennDamm11}
P.~Benner and T.~Damm, ``Lyapunov equations, energy functionals, and model
  order reduction of bilinear and stochastic systems,'' \emph{{SIAM} J. Control
  Optim.}, vol.~49, no.~2, pp. 686--711, 2011.

\bibitem{BennDamm17}
P.~Benner, T.~Damm, and Y.~{Rocio Rodriguez Cruz}, ``Dual pairs of generalized
  {L}yapunov inequalities and balanced truncation of stochastic linear
  systems,'' \emph{IEEE Transactions on Automatic Control}, vol.~62, no.~2, pp.
  782--791, 2017.

\end{thebibliography}

\end{document}